\documentclass[a4paper]{article}

\usepackage{amsmath,amssymb,amsfonts,amsthm,amscd,makeidx,stackrel,stmaryrd,float,epsfig,multicol}
\usepackage[left=3cm,top=2.5cm,bottom=2.5cm,right=3cm]{geometry}
\usepackage{amsthm}
\usepackage{array}
\usepackage{bm}
\usepackage{eucal}
\numberwithin {equation}
{section}
\theoremstyle{definition}
\newcounter{dummy} \numberwithin{dummy}{section}

\newtheorem{osservazione}[dummy]{Remark}

\numberwithin{equation}{section}
\usepackage{graphicx}
\title{Hypergeometric solutions to a three dimensional dissipative oscillator driven by aperiodic forces}
\author{Alessio Bocci, Giovanni Mingari Scarpello, Daniele Ritelli}
\date{}

\begin{document}
\maketitle

\begin{abstract}
We model the dynamical behavior of a three dimensional  (3-D) dissipative oscillator consisting of a $m$-block whose vertical fall occurs against a spring and which can also slide horizontally on a rigid truss rotating at an assigned angular speed $\omega(t)$. The bead's $z$-vertical time law is obvious, whilst its $x$-motion along the horizontal arm is ruled by a linear  differential equation we solve through the Hermite functions and the Kummer \cite{kum} confluent Hypergeometric Function (CHF)  $_{1}F_{1}$ .  After the rotation $\theta(t)$ has been computed, we  know completely the $m$-motion in a cylindrical frame reference so that some transients have then been analyzed. Finally, further effects as an inclined slide and a contact dry friction have been added to the problem, so that the motion differential equation becomes inhomogeneous: we resort to Lagrange method of  variation of constants, helped by a Fourier-Bessel expansion, in order to manage the relevant intractable integrations.


\end{abstract}
\section{Introduction}
\subsection{Aim of the problem}

Almost all types of oscillations have been the subject of intensive research for many years and several methods have been used to find exact/approximate solutions to such dynamical systems.
The oscillator hereinafter analyzed, is moving along three directions: its motion (main) ordinary differential equations is linear with variable coefficients.
Such a dissipative oscillator undertakes: own weight, constraint reaction, Coriolis force, centrifugal force, elastic forces, viscous resistance and the inertia force due to the angular acceleration $\dot{\boldsymbol{\omega}}$. All the analytical evaluations have been performed without any approximation. After that dry friction has been added, the particular integral of the consequent inhomogeneous equation could only be found by means of  variation of constants. Let us introduce the system main features.

A rotating blade is modeled by a rod of length $2L$ on which a mass $m$ is sliding, secured to the hub $\Omega$ by a spring of stiffness $k_{2}$.
The rod is jointed at $\Omega$ to a  vertical shaft of length H: the heavy block's vertical motion is against a spring of stiffness $k_{1}$.
Therefore the bead P modeling the $m$-block, has a 3-D motion\footnote{Another paper concerning a rod 3-D problem but involving elliptic functions is \cite{Laur1}}: its position is set out by a cylindrical coordinate frame with $x, z,\theta$ depending on time. The first two co-ordinates, unknown functions of time, are the object of our inquiry; the third one is decided by $\omega(t)$.
Such a model, with an obvious change, will also describe the effects induced by an increasing angular speed $\omega$ from a rest up to an absolute maximum. 

We put a reference frame $(x,\,y,\,z)$ co-mobile with the rod, so that the $x-$motion of the block, which can move on both sides with respect to the hub $\Omega$, is referred to such a non-inertial frame. In such a case, the Coriolis theorem here arises in its complete six terms formulation:
\begin{equation} \label{Coriolis}
\textbf{a}_r=\textbf{a}_a-\textbf{a}_\Omega- \boldsymbol{\omega}\wedge( \boldsymbol{\omega} \wedge\     \boldsymbol{r})-2\boldsymbol{\omega}\wedge\boldsymbol{v}_r- \dot{\boldsymbol{\omega}}\wedge \boldsymbol{r}.
\end{equation}
 The problems whose closed form solution comes out taking {\it all the terms} of Coriolis theorem are presumably very few: we cannot cite any one of them.
  
The CHFs are usually occurring in Mathematical Physics after splitting a linear partial differential equations in more ordinary differential equations. Their use for  linear or nonlinear Mechanics problems is not much practiced. We can refer for instance to \cite{Laur4}, \cite{Cit}, \cite{la}. On the contrary, much greater is their appearance in quantum mechanics: in \cite{sea} we meet CHFs in treatments on the one-dimensional harmonic oscillator and in the isotropic one. Monographic, quite old, but authoritative books are \cite{sla} and \cite{tric} as well. 
An extensive more recent introduction to CHFs can be found in \cite{ge} which holds a bibliography of 89 entries.
In \cite{Nagar}, several properties of extended CHFs are studied showing how they are met in statistical distribution theory.

\subsection{The mathematical model}
\begin{figure}[h]
\begin{center}
\scalebox{0.4}{\includegraphics{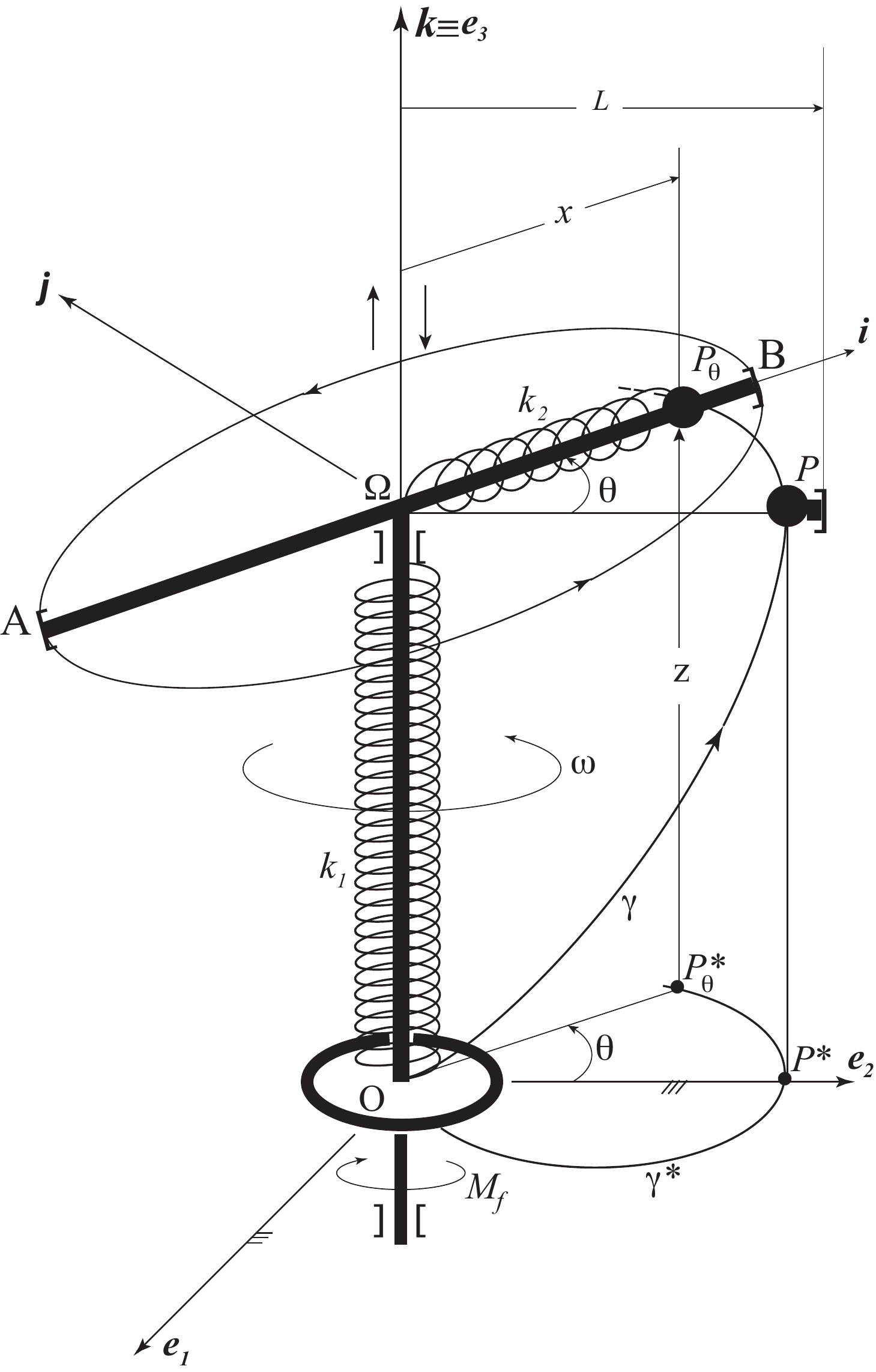}}
\caption{A 3-D oscillator and its co-mobile reference frame.}\label{f01}
\end{center}
\end{figure}
With reference to \figurename~\ref{f01}, let the mobile frame be rotating around a fixed one: if  $\boldsymbol{k}$ is the unit vector on the vertical direction $\Omega z$, the angular speed will be
$
\boldsymbol{\omega} =\omega_0 (1-qt) \boldsymbol{k}
$ 
so that the angular acceleration is: $\dot{\boldsymbol{\omega}}=-\omega_0 q\boldsymbol{k}$. 
The third co-mobile is the $\Omega y$-axis of unit vector $\boldsymbol{j}$ orthogonal to the $\Omega x$ one.
Multiplying to $m$, mass of the bead P, by \eqref{Coriolis} we get the forces a co-mobile observer sees acting on the bead. They are:

\vspace{3mm}
 
$\bullet$ the force due to the motion of the mobile reference system with respect to the fixed one: $ m\boldsymbol{a}_\Omega=m \ddot{z}\boldsymbol{k},$
  
$\bullet$ the \lq\lq absolute'' force given by $m\boldsymbol{a}_a =(-k_2x-A \dot{x}(t))\boldsymbol{i}+R_y\boldsymbol{j}+(R_z-mg) \boldsymbol{k},$ being the resistance $A \dot{x}(t)$ produced by a viscous damper, with $A>0$ drag coefficient. Notice that $\boldsymbol {R}=(0,R_{y}, R_{z})$ is the constraint unknown reaction from the rod on the $m$-bead,

$\bullet$ the centrifugal force: $m \boldsymbol{\omega}\wedge(\boldsymbol{\omega}\wedge\ \boldsymbol{r}) = -m\omega_0^2  (1-qt)^2 x \boldsymbol{i},$
 
$\bullet$ the Coriolis force: $m\boldsymbol{\omega}\wedge\boldsymbol{v}_r = m \omega_0 (1-qt) \dot{x} \boldsymbol{k},$

$\bullet$ eventually, the inertia force due to the angular acceleration:
 $m\dot{\boldsymbol{\omega}} \wedge \boldsymbol{r} =-m \omega_0qx \boldsymbol{k}. $ 

\vspace{3mm}

By \eqref{Coriolis}, we get:
\begin{equation}\label{equadin}
m\begin{bmatrix}
\ddot x\\0\\\ddot z 
\end{bmatrix}=-m\begin{bmatrix}
-\omega_0^2  (1-qt)^2 x\\ 2\omega_0 (1-qt) \dot{x}-\omega_0qx \\ 0
\end{bmatrix} +\begin{bmatrix}
-k_2x\\ R_y\\ R_z-mg
\end{bmatrix} +\begin{bmatrix}
A \dot{x}(t)\\ 0\\ 0
\end{bmatrix} 
\end{equation}
namely:
\begin{equation}\label{equadinn}
\begin{cases} \ddot x+A \dot{x}(t) +\Biggl[\dfrac{k_2}{m}-\omega_0^2  (1-qt)^2\Biggl] x=0 \\ R_y = 2m\omega_0 (1-qt) \dot{x}-m\omega_0qx \\
R_z = mg +m \ddot{z}
 \end{cases}
\end{equation}

\begin{osservazione}
After adding the second spring,  \eqref{equadinn} shows that as more the origin of the co-mobile frame accelerates with respect to the origin of the fixed one, the more increases the reaction component $R_z$. 

\end{osservazione}
\begin{osservazione}
The first of \eqref{equadinn} provides a movement pure equation; the second one gives the $y$ component of constraint reaction after the motion has been solved; the third will allow to evaluate the $z$ component of constraint reaction, after finding the time law $z=z(t)$  of the arm $\Omega P$ oscillations.  
\end{osservazione} 

\section{From $\boldsymbol{x}$-the motion equation to the Hypergeometric Confluent function}
Let us integrate the first of \eqref{equadinn}, for shortness represented as:
\begin{equation}\label{ehh}
\ddot{x}+A \dot{x}-(at^2+bt+c)x=0 \\
\end{equation}
where $a=q^2 \omega_0^2 >0,\,A>0,\, b=-2q\omega_0^2<0$ and $c=\omega_0^2-\tfrac{k_{2}}{m}.$ We state that the general solution of \eqref{ehh}, putting $\beta=\tfrac{b^2-a(A^2+4c)}{8a^{3/2}}$ is given by the formula:
\begin{equation}
x(t)=e^{-\frac{at^2+t(b+\sqrt{a}A)}{2\sqrt{a}}}\left[C_1\,H_{\beta-\frac{1}{2}}\left(\frac{b+2at}{2a^{3/4}}\right)+C_2\,_1{\rm F}_1\left( \left. 
\begin{array}{c}
\frac{1}{4}-\frac{\beta}{2} \\[2mm]
\frac{1}{2}
\end{array}
\right| \frac{(b+2at)^2}{4a^{3/2}}\right)\right]
\label{eh}
\end{equation}
where $H_{\alpha}(t)$ is the Hermite function of nonintegral order $\alpha$ and argument $t,$ while $_1{\rm F}_1$ is the Kummer function; see for both the Appendix at the end of the paper for some details.

To obtain the analytic solution \eqref{eh} we make twice a variable transformation, the first involving the dependent variable and the second concerning the independent variable, leading to an Hermite differential equation of non integer order for the unknown function $\sigma=\sigma(\tau)$ see formula \eqref{hermy} in the Appendix. The first transformation of equation \eqref{ehh} is obtained putting 
\begin{equation}\label{changeof}
x(t)=s(t) \exp\left({-\frac{\sqrt{a} A t+a t^2+b t}{2 \sqrt{a}}}\right)
\end{equation}
so that \eqref{ehh} is changed in:
\begin{equation}\label{quaqua}
\ddot{s}-\frac{(2 a
   t+b) }{\sqrt{a}}\dot{s}-\frac{\left(4 a^{3/2}+a A^2+4
   a c-b^2\right)}{4 a}s=0
\end{equation}
getting rid of the second degree term in $t$. Thus \eqref{quaqua} is of the following type:
\begin{equation}\label{quaquagen}
\ddot{s}-(2m_2 t+m_1)\dot{s}+2m_0 s=0.
\end{equation}
where we put:
\begin{equation}\label{parasta}
m_2=\sqrt{a},\,m_1=\frac{b}{\sqrt{a}},\,m_0=\frac{\left(b^2-4 a^{3/2}-a A^2-4
   a c\right)}{8 a}
\end{equation}
The second step will drive \eqref{quaquagen} into an Hermite equation, see \eqref{hermy} in the Appendix, by putting
\[
t=\frac{\tau}{\sqrt{m_2}}-\frac{m_1}{2m_2}\iff \tau=\sqrt{m_2}t+\frac{m_1}{2\sqrt{m_2}}
\]
we obtain the Hermite differential equation, see again \eqref{hermy} in the Appendix, where we write $\sigma(\tau)$ for $s\left(\frac{\tau}{\sqrt{m_2}}-\frac{m_1}{2m_2}\right)$:
\begin{equation}\label{quaquagen2}
\sigma''-2\tau\sigma'+2\frac{m_0}{m_2}\sigma=0
\end{equation}
Solution \eqref{eh} comes then form formula \eqref{hermysol} combined with \eqref{changeof} and \eqref{parasta}.

The analytical solution of equation \eqref{ehh} of our model is then given by \eqref {eh}, where:
\begin{equation}\label{param}
 a=\omega_0^2q^2>0,\, b=-2q\omega_0^2<0,\, c=\omega_0^2-\frac{k_2}{m}\lesseqgtr0. 
 \end{equation}
Writing the Lagrange $z$-motion equation, we get the initial value problem:
\begin{equation}
\begin{cases} 
\ddot{z}+\dfrac{k_1}{m}z=-g \\[2mm] 
z(0)=z_0 \\ 
\dot{z}(0)=\dot{z}_0
\end{cases} 
\end{equation}
which is solved as:
\begin{equation}\label{zedt}
z(t)=C_3\sin\left(\overline{\omega}\,t+C_4\right)-\frac{mg}{k_{1}}
\end{equation}
with $\overline{\omega}=(k_1/m)^{1/2}$. By the initial conditions $z(0)=z_0,\, \dot{z}(0)=\dot{z}_0,$ we get:
\[
\begin{cases} 
C_3=\sqrt{\Biggl(z_0+\dfrac{mg}{k_1}\Biggl)^2+\dfrac{\dot{z}_0^2}{\overline{\omega}^2}} \\ 
C_4=\arctan\Biggl(\dfrac{\overline{\omega}(k_1z_0+mg)}{\dot{z}_0k_1}\Biggl) 
\end{cases} 
\]
Inserting $z(t)$ given by \eqref{zedt} into the third of \eqref{equadinn}, we get the $z$-component of the reaction:
\begin{equation}
R_z=m[g-C_3\overline{\omega}^2 \cos(\varpi t+C_4)].
\end{equation}
The $y$-component of the reaction requires, see the second of \eqref{equadinn}, $x$ and its time derivative:
\begin{equation}
R_y = 2m\omega_0 (1-qt) \dot{x}(t)-m\omega_0qx(t).
\end{equation}
 The Kummer derivative theorem provides such a derivative, but we omit to relate here its too long expression.

To complete our description, we point out how the polar anomaly $\theta$ can be expressed as a time-function.
The angular speed is changing during time: the rod $O\Omega$ is first pinning with angular speed $\omega_0$ which, starting from $t=0$ is affected by a  braking torque $ M_f$:
\begin{equation}
\begin{cases} J \dot \omega=-M_f \\ \omega(0) =\omega_0\end{cases} 
\end{equation}
So that the angle $\theta$ changes in time as
\begin{equation}
\begin{cases} \dfrac{{\rm d}\theta}{{\rm d}t}=\omega_{0}(1-qt) \\[2mm] 
\theta(0) =0\end{cases} 
\end{equation}
where $q=M_f/J_z\omega_{0}$ being $J_z$ the inertia axial moment of the shaft.

The reversion of the polar anomaly $\theta=\theta(t)$ provides time as a function of the instantaneous angle, $t=t(\theta)$:
\begin{equation}\label{t-theta}
\begin{cases} t(\theta)=\dfrac{1}{q}\left(1-\sqrt{1-\dfrac{2q\theta}{ \omega_0}}\right), \\[4mm] 0\leq\theta\leq  \dfrac{\omega_0}{2q} \end{cases} 
\end{equation}

The planar line $\gamma^*$ is the projection on the fixed plane ($\boldsymbol{e_1},O, \boldsymbol{e_2})$, see \figurename~\ref{f04} of the $3-D$ torse $\gamma$. Its equation follows by inserting $t=t(\theta)$ from \eqref{t-theta} within $x(t)$ which, because of the changed framework, will be renamed as $\rho$:
\begin{equation}\label{t-theta2}
\begin{cases} \rho=\rho(\theta), \\[2mm]
0\leq\theta\leq  \dfrac{\omega_0}{2q}. 
\end{cases} 
\end{equation}
So that the polar equation of  $\gamma^*$  is obtained.

\section{A sample problem}

The $x$-problem we faced, is a initial value one with two initial conditions $x_0$ and $v_0$, so that the integration constants $c_1$ and $c_2$ will depend on: $a, b, c, x_0$ and $v_0$. Then the $x$-solution to \eqref{eh} will be of the kind 
$
x=f (a, b, c; x_0, v_0; t).
$
Being the Kummer function unbounded, and in lack of any barrier to $x$, in carrying out the simulations, it could be possible to see some $x$ behavior not fit with its boundaries.

First of all then we have to state a time span since the forcing transient expires at $t=t_{\lim}=1/q$. 
Afterwards, we should subject the function $x(t)$ of \eqref{eh} to a constraint so that $x$ could go never outside the range $(-L, L)$.
Notice that $x$ is provided by the \eqref{eh} depending on six variables: then the algorithmic complexity of such a formulation and its low practical usefulness advised us against improving it, so that we restrict here to show how the boundaries could be acted:

\begin{equation}
\begin{cases} -L\leq x (a, b, c; x_0, v_0; t)\leq L \\ 0\leq t\leq 1/q\end{cases} 
\end{equation}
In practice, we assumed as a sample problem, the following physical data, stopping the computation whenever should $x(t)$ go out its boundaries.
\begin{multicols}{3}
\begin{description}
\item $L=1m; $
\item $H=3m;$
\item $m=1 kg;$
\item $k_1=10Nm^{-1}$;
\item $k_2=Nm^{-1},\, {\rm variable}$; 
\item $\omega_0=3 {\rm rad}s^{-1};$
\item $q=s^{-1},\, {\rm variable};$
\item[]
\end{description}
\end{multicols}
\noindent where the $k_2$ values ($Nm^{-1}$) for transients I, II, III, IV, V respectively, are: 10; 8; 30; 8, .... The $q$ values ($s^{-1}$) for transients I, II, III, IV respectively, are: 1/10, 1/10, -1/10,-1/10. We carried out five transients.

Notice that $c=\omega_0^2-\tfrac{k_2}{m}$ has dimension of a square of a frequency and sets up a comparison on the resulting bead motion between the centrifugal outside cause $\omega_0$ and the spring opposition measured by $\tfrac{k_2}{m}$. Therefore the most studied occurrence $c<0$ (transients I, II, III, V) means that the elastic influence is prevailing on the rotational starting speed.

\subsection{I transient: $\omega$ shutdown with major elastic force influence}

\begin{figure}[h]
\begin{center}
\scalebox{0.35}{\includegraphics{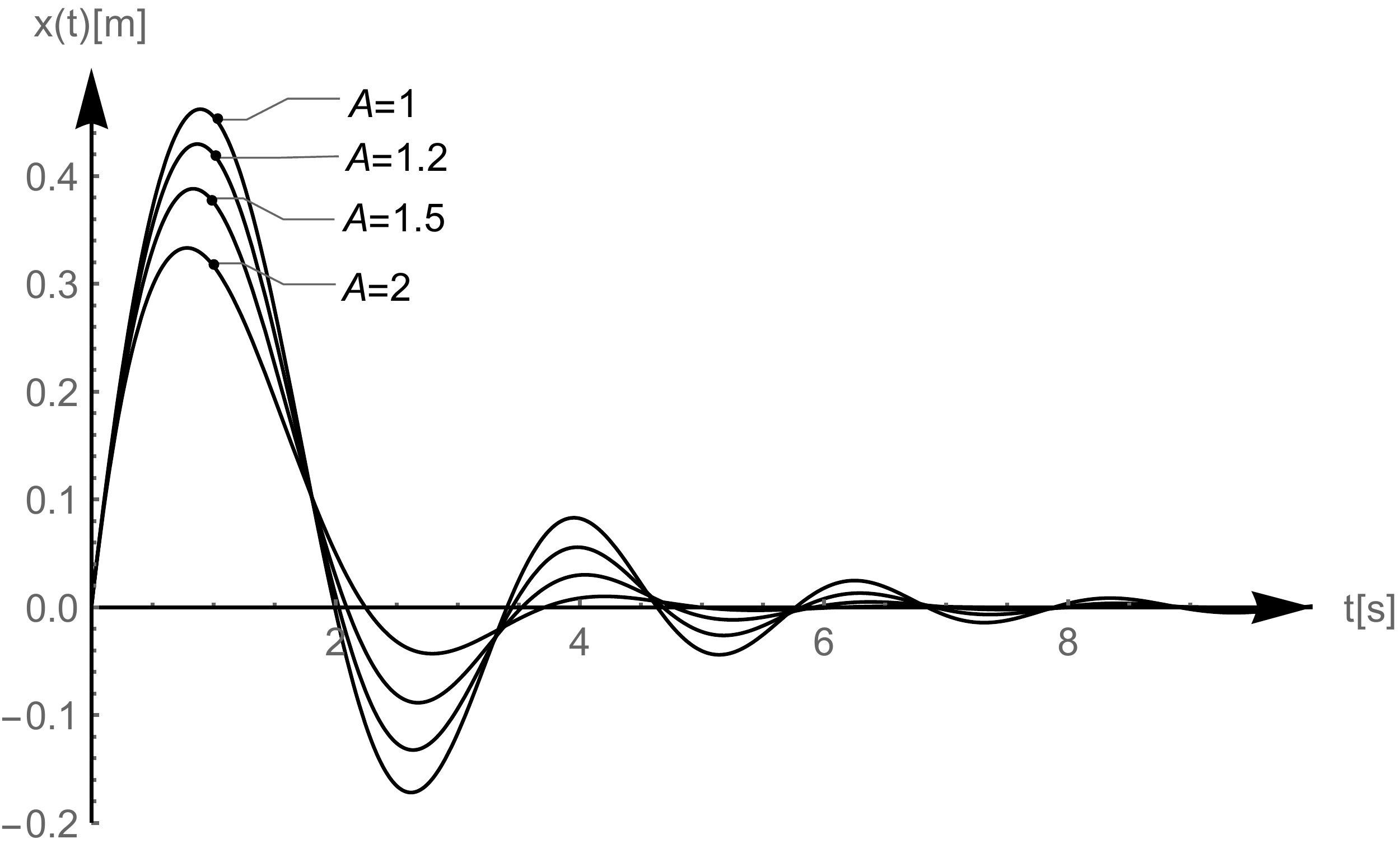}}
\caption{I transient: $\omega$ shutdown with prevailing elastic force ($c<0$): $(x, t)$-oscillations under four different A drag values.}\label{f02}
\end{center}
\end{figure}

By formulae \eqref{param} we see that both $a$ and $b$ signs do not depend on $k_2$, while the mix of values of $m,  k_2, \omega_0$ decides that of $c=\omega_0^2-\frac{k_2}{m}$.
Let be $\frac{k_2}{m}>\omega_0^2$: then the elastic spring force will be, at the beginning, of greater effect than the centrifugal one, so that $x$-oscillations will occur and will expire with the angular speed, see \figurename~\ref{f02}

\begin{figure}[H]
\begin{center}
\scalebox{0.19}{\includegraphics{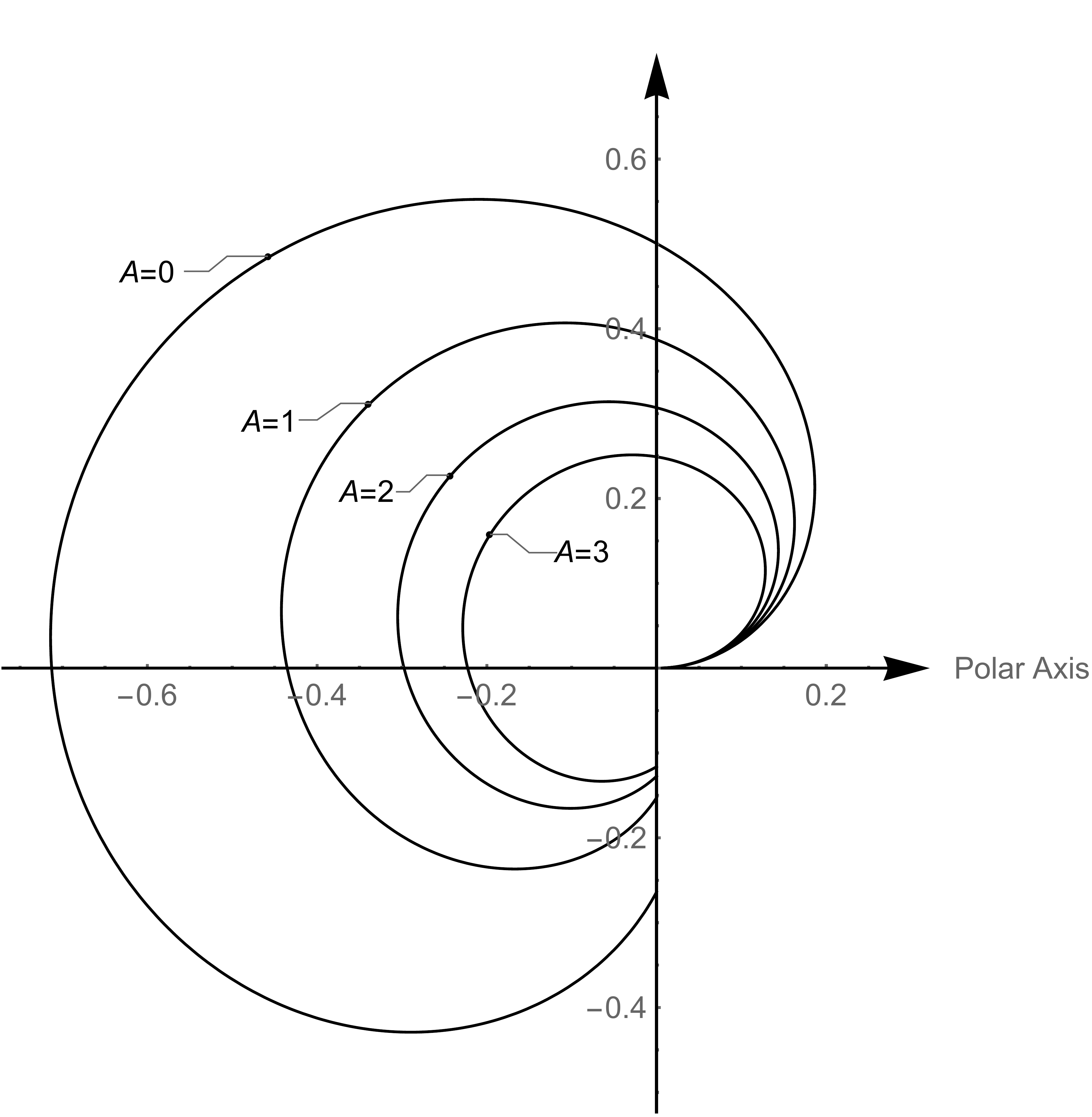}}
\caption{$\omega$ shutdown with prevailing elastic influence: polar plot $\rho=\rho(\theta)$ of the $\gamma^{\star}$ projection of the bead's torse. Four different $A$-drag values.}\label{f04}
\end{center}
\end{figure}

We can also compute by the second of \eqref{equadinn} the oscillations of reaction $R_y$ which is affecting the bead:  see \figurename~\ref{f03}. Putting $t= t (\theta)$ within $x(t)$, we get the polar plot, to be read counterclockwise, of the curve $\gamma^*$ projection of the torse $\gamma$ on the fixed plane ($\boldsymbol{e_1},O, \boldsymbol{e_2})$, see \figurename~\ref{f01} and \figurename~\ref{f04}.

\subsection{II transient: $\omega$ shutdown with minor elastic influence: oscillations.}
\begin{figure}[H]
\begin{center}
\scalebox{0.45}{\includegraphics{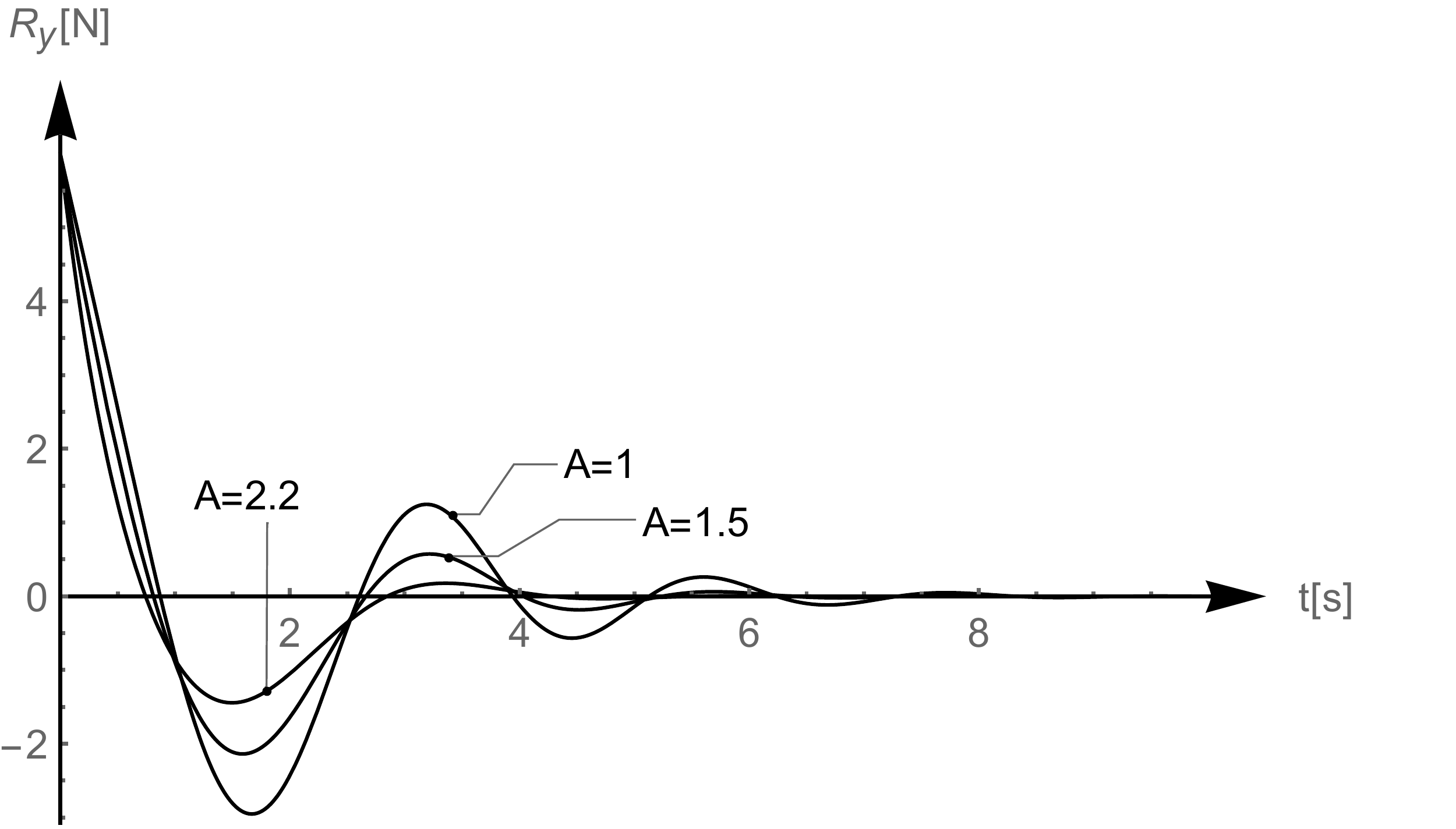}}
\caption{II transient: $\omega$ shutdown with prevailing elastic force, $c<0:\,(R_y, t)$-oscillations under different $A$ drag values.}\label{f03}
\end{center}
\end{figure}

The simulations with smaller $k_2$ values produced not-oscillatory behaviors, here omitted for shortness.
Increasing  $k_2$ values, the effect of the elastic forces is not too weak, so that reactions will oscillate. We provide only the plot of $R_y$ versus $t$, \figurename~\ref{f03}
\subsection{III transient: $\omega$ increasing with major elastic influence}
\begin{figure}[h]
\begin{center}
\scalebox{0.25}{\includegraphics{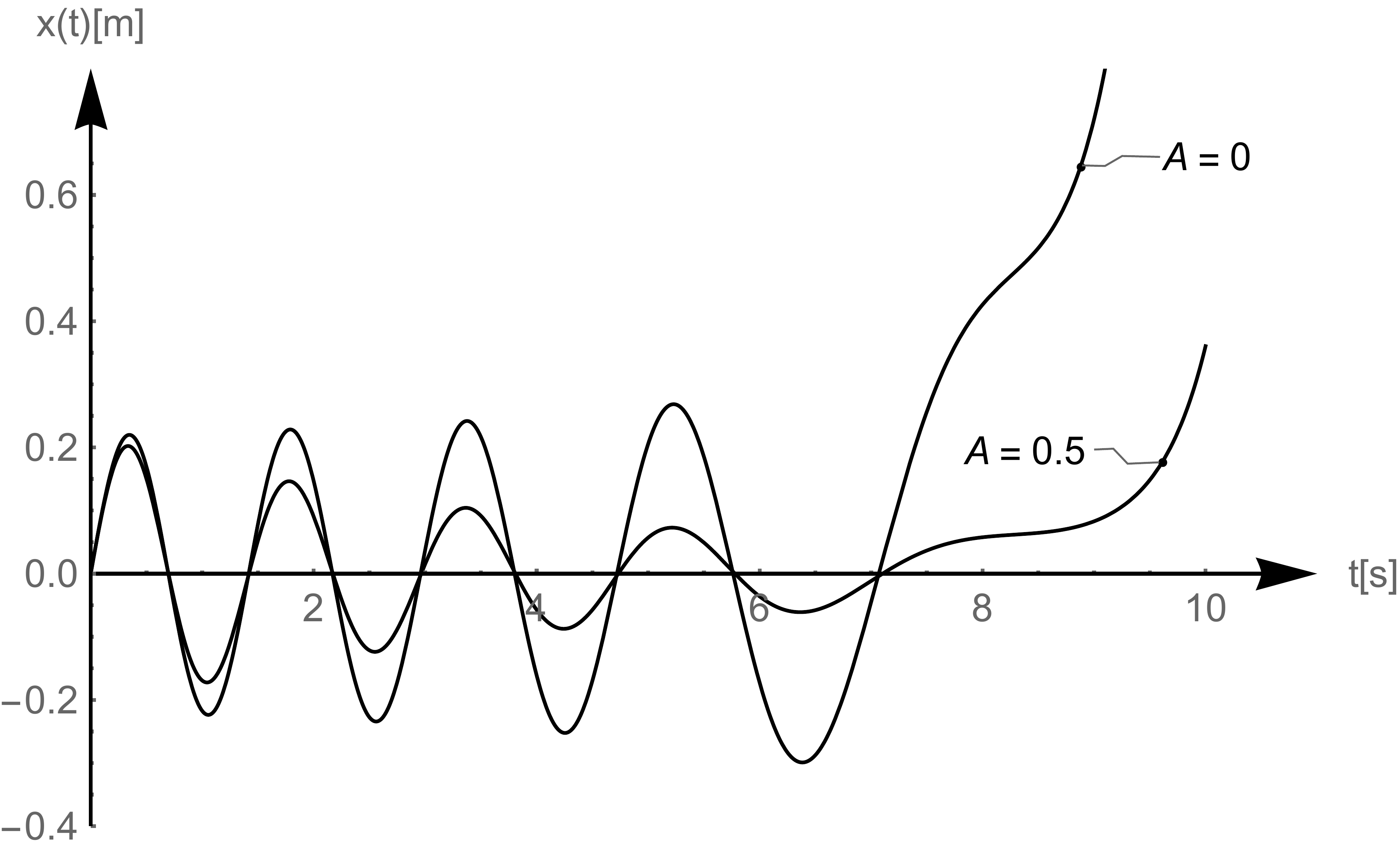}}
\caption{III transient: $\omega$ increasing with major elastic influence, $c<0:$  bead's $(x,t)$-oscillatory first course and final blow-up. Different A drag values.}\label{f06}
\end{center}
\end{figure}

By \eqref{param} we see that the $\omega$ increase, namely $q<0$, impacts on the sign of $b$, and then on the exponential factor responsible of the drop of $x$, according to \eqref{eh}. We see the major elastic influence, $c<0,$  succeeds in producing decaying oscillations again, see \figurename~\ref{f06}, but at the end the $(x, t)$ plot will blow-up in any case.

\subsection{IV transient: $\omega$ increasing with minor elastic influence, non-oscillatory bead's $x$-behaviour}

\begin{figure}[h]
\begin{center}
\scalebox{0.25}{\includegraphics{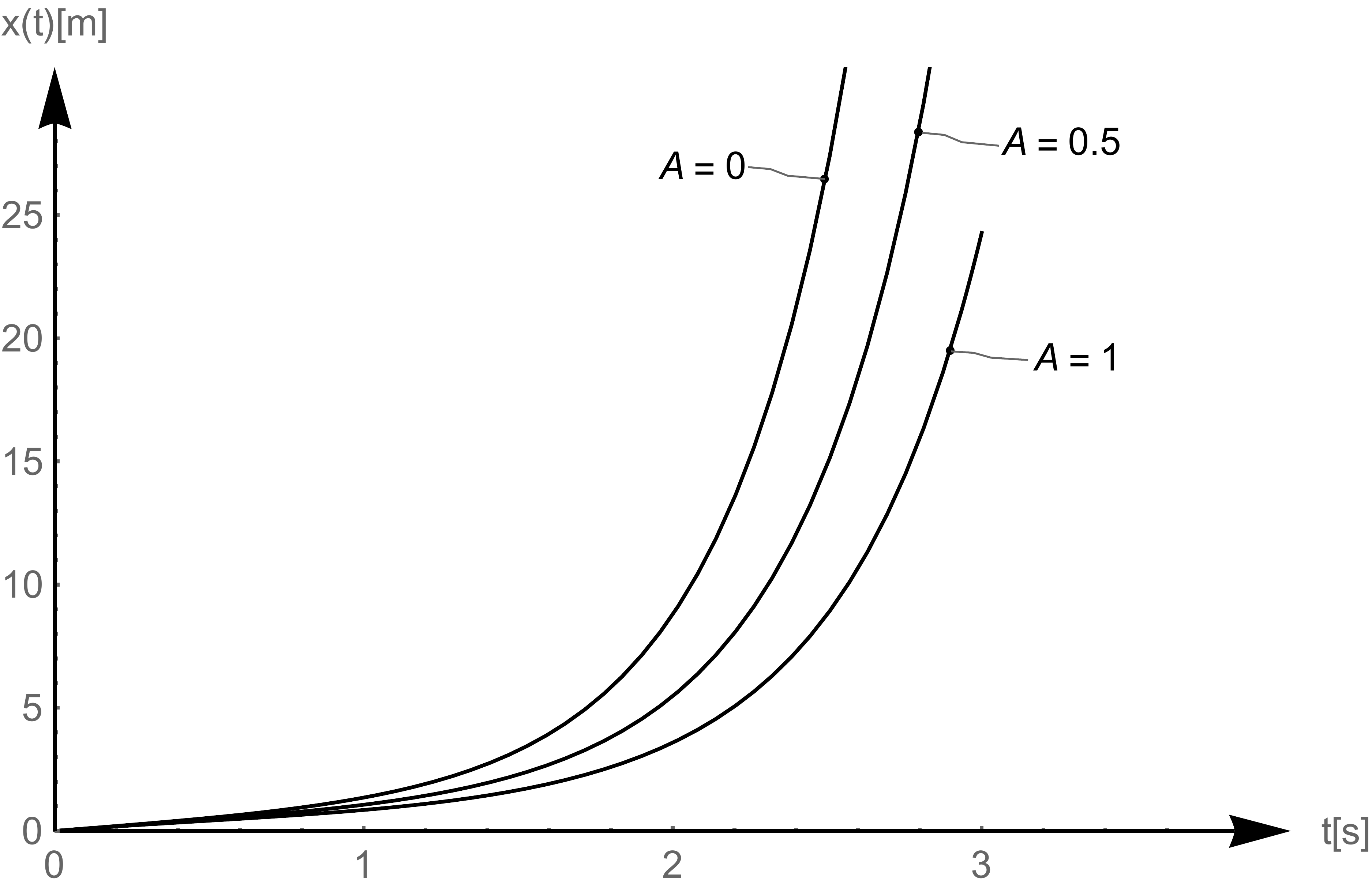}}
\caption{IV transient: $\omega$ increasing with not-oscillatory bead's $(x, t)$-plot. Different A drag values.}\label{f07}
\end{center}
\end{figure}

Here, see \figurename~\ref{f07}, we can see how the double event of the $\omega$ growth and the minor elastic influence, $c>0,$ will produce an $x$-monotonic path whose plot has been cut-off about the third second due to its blowing-up.

We know that oscillations occur when in a system's structure the elastic features are the prevailing ones. \figurename~\ref{f06} shows how this happens for the first times when $\omega$ is not yet so large. Next, when the rotational effects will prevail, the vibrations do cease at all and the motions tend to instability with or without damping which can affect only the blow-up rate.

\begin{osservazione}
Specially for those transients where $\omega,$ and then the centrifugal force, grows, the knowledge of the time law of the transverse $R_y$ of the rod-blade on the block, will provide, by reaction, the loads acted by the block on the blade itself. In such a way, a comparative check on the $\omega$ range, will define the dynamic safety levels for the rod during the start-up transients.
\end{osservazione}

\subsection{V transient}

Here we are analyzing a transient where we assume the angular speed to be a constant, $q=0,\,\omega(t)=\omega_0,$ and with $c<0$ as before, namely the elastic forcing keeps its control over the centrifugal one. We can then  expect pure $x$-oscillations whose expiring time is ruled by the damping $A$ amount. The motion equation collapses in a constant coefficients one.

\begin{figure}[h]
\begin{center}
\scalebox{0.25}{\includegraphics{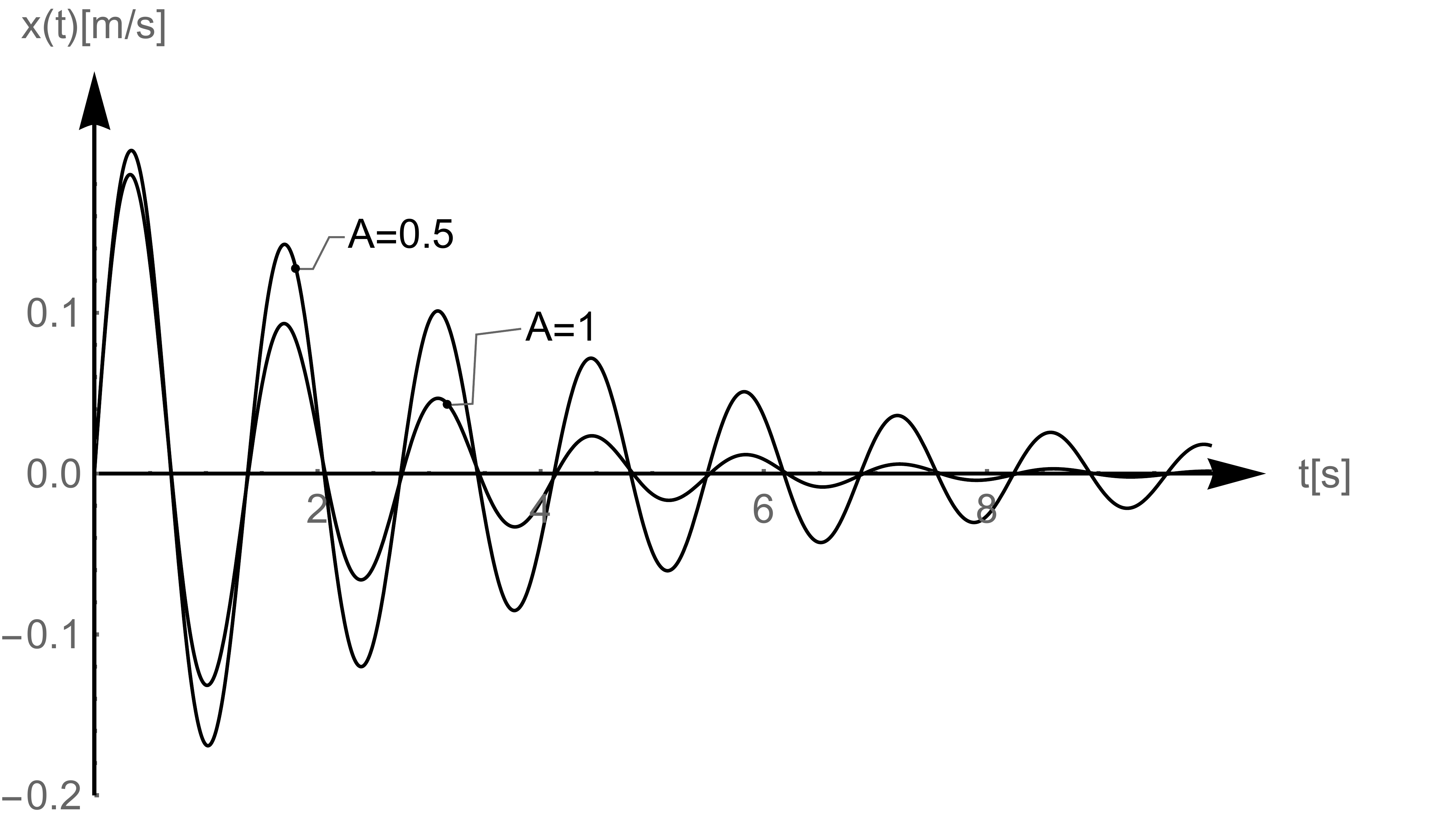}}
\caption{V transient: $\omega$ constant with mayor elastic influence, $c<0:$ damped oscillations.}\label{f0V}
\end{center}
\end{figure}

\subsection{Summary}

The main features of all transients ruled by the homogeneous Weber equation \eqref{ehh} are  summarized below

\begin{center}
    \begin{tabular}{ | c | l | c | p{5cm} | l |}
    \hline
    {\it Transient N.}&      $\omega(t)$ {\it behavior} & $c=\omega_0-\tfrac{k}{m^2}$ & {\it Main computed effects} & {\it Figure N.}\\ \hline
    I & shutdown & $c_{\rm I}<0$ & Clear $x$ oscillations, different dampings & \figurename~\ref{f02} \\ %\hline
    \cline{4-5}
  &  &  & The head's curve projections are spiralizing & \figurename~\ref{f04}\\ \hline

    II & shutdown & $c_{{\rm II}}<c_{\rm I}<0$ & Oscillations of constraint reaction $R_y$ & \figurename~\ref{f03} \\ \hline
    III & increasing & $c_{{\rm III}}<0$ & Initially $x$ oscillates but at end blows up & \figurename~\ref{f06} \\ \hline
    IV & increasing & $c_{{\rm IV}}>0$ & The $x$ monotonic increase plot & \figurename~\ref{f07} \\ \hline
    V & constant $=\omega_0$ & $c_{{\rm V}}<0$ & Some $x$ oscillations without blow-up, all expiring sooner or later according to the damping $A$ & \figurename~\ref{f0V} \\
    \hline

    \end{tabular}
\end{center}

\section{The forced case adding dry friction too}

By adding into the system an outside forcing $\mu$ which models both effects of a blade inclination and a dry friction contact, the previous motion equation becomes inhomogeneous:
\begin{equation}
\begin{cases}
\ddot{x}(t)+A \dot{x}(t)-(at^2+bt+c)x(t)=\mu\\
x(0)=0\\
\dot{x}(0)=1,
\end{cases}
\end{equation} 
with $\mu \in \mathbb{R}$ is homogeneous to an acceleration.
Applying the Lagrange variation of parameters method\footnote{Between 1778 and 1783, Lagrange developed his method  in papers concerning the variations in planetary motions and in another series of memoirs on computing the orbit of a comet from three observations. During 1808-1810, he gave to the method its final form in a series of papers among which we shall quote \cite{lagrangia} which is available at Gallica web site}, the particular integral will be provided by:
\begin{equation}
\overline{x}(t)=c_1(t)x_1(t)+c_2(t)x_2(t)
\label{lola}
\end{equation}
being $c_1(t)$ and $c_2(t)$ given by:
\begin{equation}
c_1(t)= -\mu\int \frac{x_2(t)}{\dot{x}_2(t)x_1(t)-x_2(t)\dot{x}_1(t)}{\rm d}t,
\label{primo}
\end{equation}
and:
\begin{equation}
c_2(t)=\mu\int\frac{x_1(t)}{\dot{x}_2(t)x_1(t)-x_2(t)\dot{x}_1(t)}{\rm d}t.
\label{secondo}
\end{equation}
Mind that $x_1(t)$ and $x_2(t)$ are two linearly independent integrals of the homogeneous equation which can be read on (\ref{eh}):
\begin{equation}
x_1(t)=e^{-\frac{at^2+t(b+\sqrt{a}A)}{2\sqrt{a}}}H_{\beta-\frac{1}{2}}\left(\frac{b+2at}{2a^{3/4}}\right)
\end{equation}
and
\begin{equation}
x_2(t)=e^{-\frac{at^2+t(b+\sqrt{a}A)}{2\sqrt{a}}}\,_1{\rm F}_1\left( \left. 
\begin{array}{c}
\frac{1}{4}-\frac{\beta}{2} \\[2mm]
\frac{1}{2}
\end{array}
\right|  \frac{(b+2at)^2}{4a^{3/2}}\right).
\end{equation}
In order to avoid a fully numerical solution we decided an approach by which, even if resorting to  definite numerical integrations, we are however allowed to keep the control of the solution's functional structure.
Both integrals (\ref{primo}), (\ref{secondo}) cannot be computed analytically: nevertheless we can resort to a series expansion of the integrand. Let us see first the function $c_1(t)$.

 Formula (\ref{primo}) is of kind:
\begin{equation}
c_1(t)=-\frac{8a^2 \mu}{4a^{3/2}-b^2+a(A^2+4c)} \int f(t){\rm d}t,
\label{sch}
\end{equation}
where:
\begin{equation}
f(t)=\frac{e^{\frac{at^2+t(b+\sqrt{a}A)}{2\sqrt{a}}}}{u_1(t)+u_2(t)}\,_1{\rm F}_1\left( \left. 
\begin{array}{c}
\frac{1}{4}-\frac{\beta}{2} \\[2mm]
\frac{1}{2}
\end{array}
\right| \frac{(b+2at)^2}{4a^{3/2}}\right)
\end{equation}
being in it:
\begin{equation}
u_1(t)=2a^{3/4}H_{\beta-\frac{3}{2}}\left(\frac{b+2at}{2a^{3/4}}\right)\,_1{\rm F}_1\left( \left. 
\begin{array}{c}
\frac{1}{4}-\frac{\beta}{2}\\[2mm]
\frac{1}{2}
\end{array}
\right| \frac{(b+2at)^2}{4a^{3/2}}\right)
\end{equation}
and:
\begin{equation}
u_2(t)=(b+2at)H_{\beta-\frac{1}{2}}\left(\frac{b+2at}{2a^{3/4}} \right)\,_1{\rm F}_1 \left( \left. 
\begin{array}{c}
\frac{5}{4}-\frac{\beta}{2} \\[2mm]
\frac{3}{2}
\end{array}
\right| \frac{(b+2at)^2}{4a^{3/2}}\right).
\end{equation}
We now go on by showing how $c_1(t)$ can be evaluated by means of a sample problem.
Let us choose for instance the case $A=1,\,a=9/100,\,b=-9/5,\,c=-1,\,\mu=\text{undefined}$. By (\ref{sch}) we see the solution $x(t)$ is by the physical problem defined for $t\in[0,1/q]$; so that, with $1/q=10$, the physical range representation of solution is $t\in[0,10]$. We could try to model our specific $f(t)$, say $\hat{f}(t)$, in such a range by means of a MacLaurin polynomial $T_{200} (t)$ holding 200 terms:
\begin{figure}[h]
\centering
\includegraphics[scale=0.4]{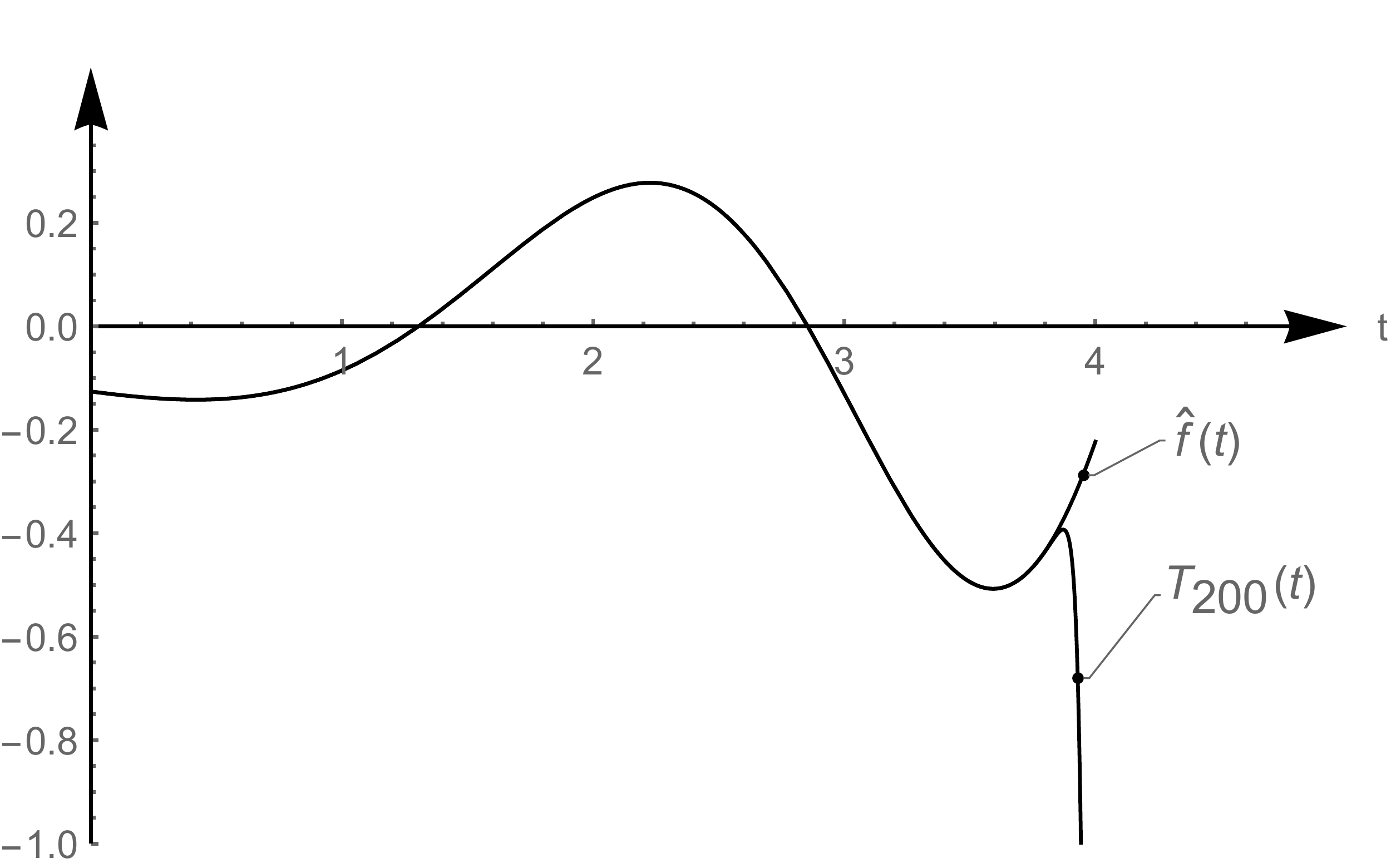}
\caption{Sketch of both curves $\hat{f}(t)$ and  $T_{200}(t)$. The ordinate computed values have been multiplied by $10^{14}$.}\label{tailor}
\end{figure}
\figurename~\ref{tailor} shows that before four seconds the approximation fails and the dissimilarity takes such a shape due to the very high powers involved in Taylor polynomial. 
We then try for $\hat{f}(t)$ a Fourier-Bessel series expansion, in terms of the Bessel function of order zero $J_0$, see \cite{watson}:
\begin{equation}
\hat{f}(t)=\sum_{k=1}^{+\infty}B_kJ_0\left(\frac{\alpha_k}{\overline{t}}t\right),
\label{be}
\end{equation}
where denoting as usual $J_0$ and $J_1$ the first kind Bessel functions of order 0 and 1 respectively. The $B_k$ coefficients can be computed as:
\begin{equation}
B_k=\frac{2}{\overline{t}^2J_1^2(\alpha_k)}\int_0^{\overline{t}}t\hat{f}(t)J_0\left(\frac{\alpha_k}{\overline{t}} t\right){\rm d}t
\label{svil}
\end{equation}
to which the condition to be added is:
\begin{equation}
J_0(\alpha_k\overline{t})=0.
\end{equation}

In previous lines $\overline{t}$ is the root of function $\hat{f}(t)$ and $\alpha_k$ is the sequence of roots of $J_0$ provided by the formula of McMahon \cite{macmahon}. We did use the function $J_0$ because it shall be $\hat{f}(0)\neq 0$.
Such a representation fits perfectly $\hat{f}(t)$ and 200 terms provide (with $A\ne 0$) a satisfactory convergence (but with $A=0$ only 30 terms are enough).
We pass to integrate $\hat{f}(t)$ as in (\ref{sch}). We choose to express the integral of $J_0$ in terms of generalized hypergeometric function $_1{\rm F}_2$% (not to be confused with the best known Gauss hypergeometric function $_2F_1$):
\begin{equation}
\int J_0(\gamma t){\rm d}t= \,_1{\rm F}_2 \left( \left. 
\begin{array}{c}
\frac{1}{2} \\[2mm]
1,\frac{3}{2}
\end{array}
\right| -\frac{\gamma^2}{4}t^2\right)
\end{equation} 
Such a relationship can be found in \cite{Luke}, page 44; inserted into the $c_1(t)$ formula, it provides:
 \begin{equation}
c_1(t)=-\frac{8a^2 \mu}{4a^{3/2}-b^2+a(A^2+4c)}\, t \, \sum_{k=1}^{+\infty}B_k \, _1{\rm F}_2\left( \left. 
\begin{array}{c}
\frac{1}{2} \\[2mm]
1,\frac{3}{2}
\end{array}
\right| -\frac{\alpha_k^2}{4\overline{t}^2} t^2 \right).
\label{ultima}
\end{equation}
So the series truncation can be done after the integration and not before, so that the approximation will result of better quality.
Let us now compute the value of $\overline{t}$, namely the root of function $\hat f(t)$ subsequent to the end ($t=10$) of the physical range of $x(t)$. We find:
$\overline{t}=10.5031$.
We shall then expect a $\hat{f}(t)$ representation converging with a relatively small number of terms. And in fact with 120 terms of the expansion \eqref{be} of  $\hat{f}(t)$ we get a satisfactory fit for $t\in[0,10]$.

Assuming $\mu=1$, by (\ref{ultima}) we get $c_1(t)$ whose plot is shown:

\begin{figure}[h]
\centering
\includegraphics[scale=0.45]{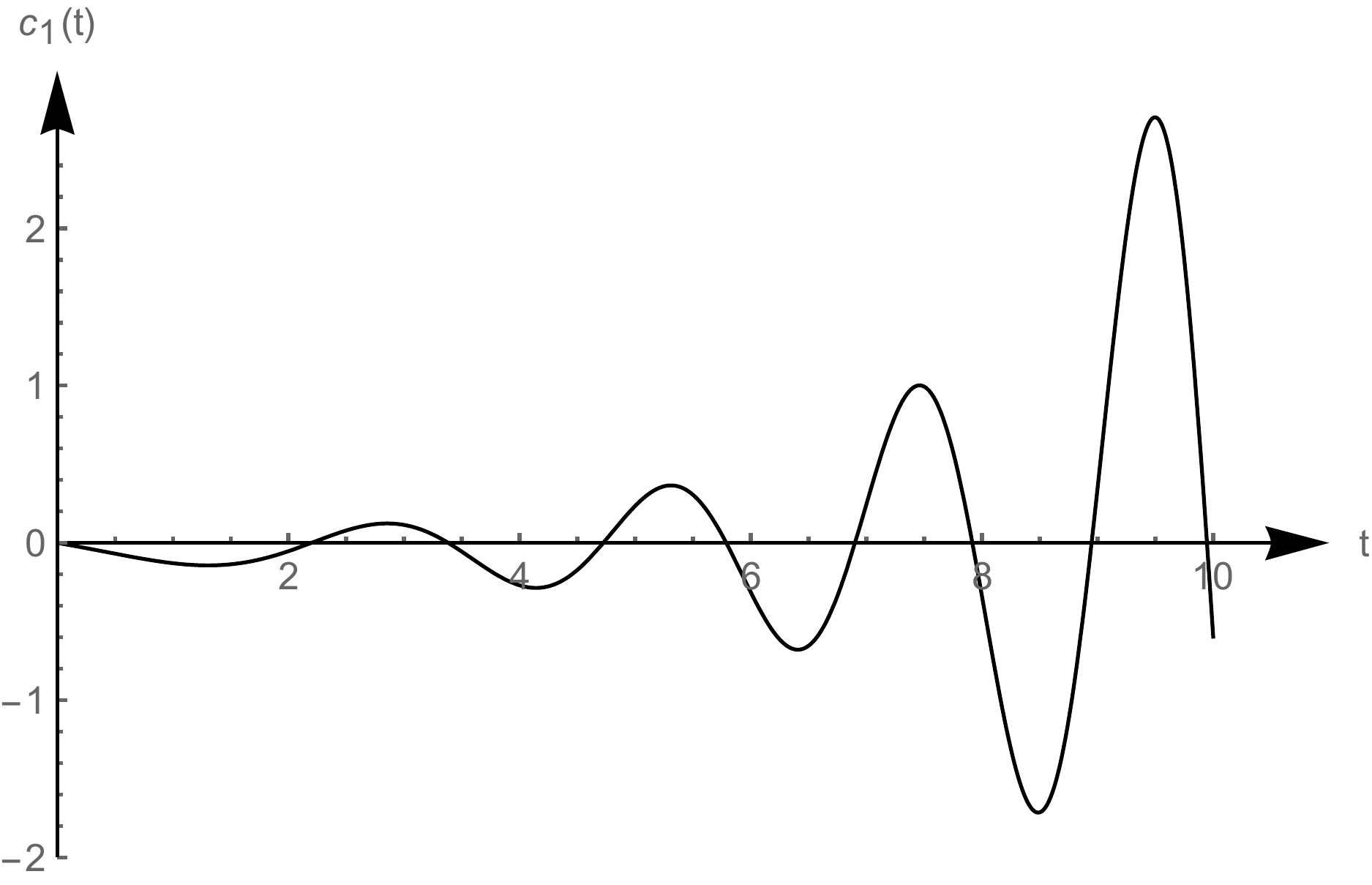}
\caption{Plot of function $c_1(t)$ to be put in \eqref{lola}: the ordinate computed values have been multiplied by $10^{14}$. }
\end{figure}

The same can be done in order to $c_2(t)$.
 The particular integral has been evaluated for $\mu=1$, but, due to formulae \eqref{primo} and \eqref{secondo}, its amplitude is just $\mu$. 
The general integral is got adding it to that of the homogeneous equation. We provide such a $\mu=1$ forced case description through some curves for different $A$ values.
\begin{figure}[h]
\centering
\includegraphics[scale=0.3]{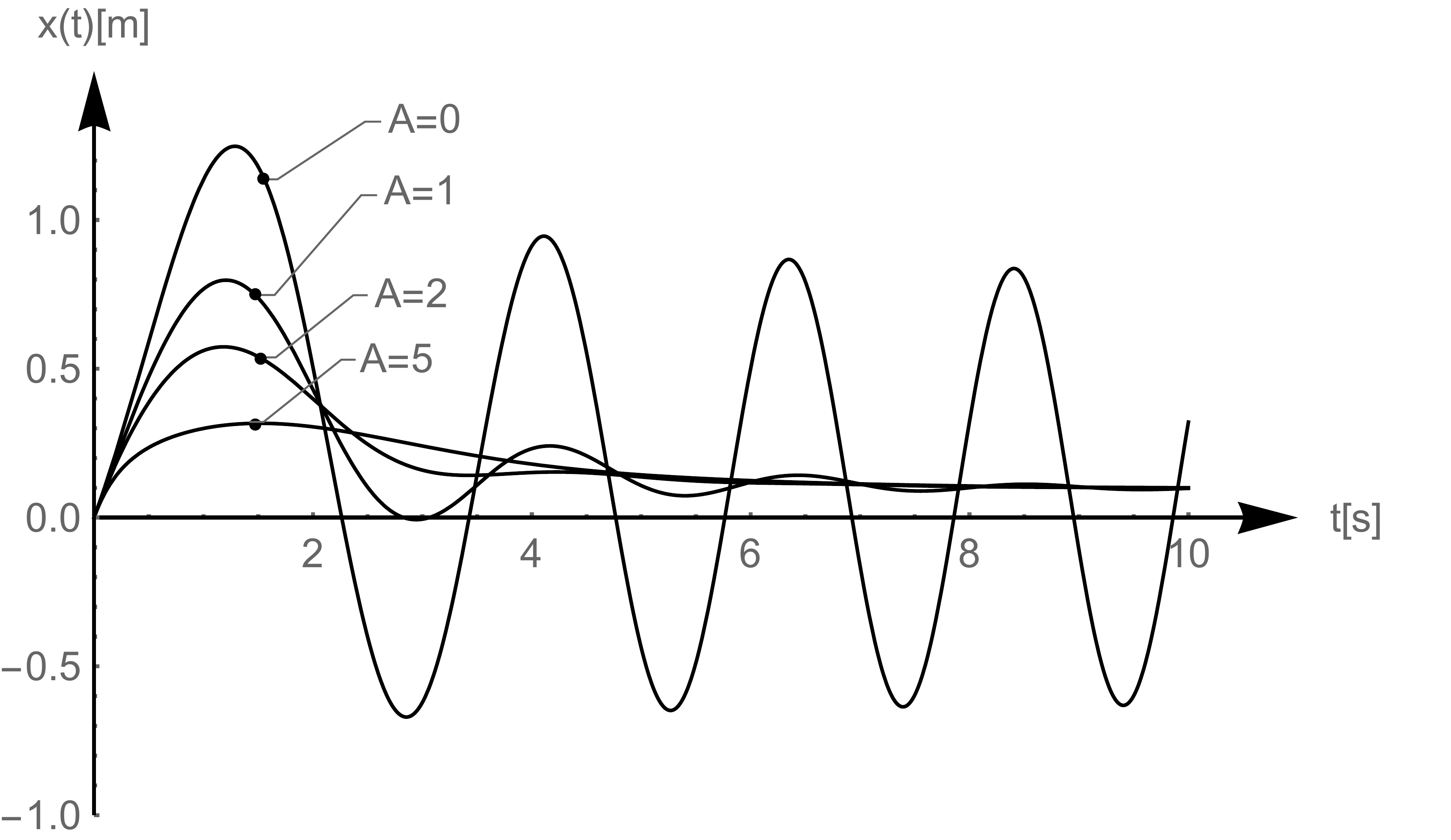}
\caption{The forced case: damped oscillations, different $A$ values and under dry friction, $\mu=1$.}\label{finove}
\end{figure}
Notice that for drag values $A> 2$ no oscillatory behavior occurs any more and  the overdamping takes place. 
\section{Conclusions}
The main features of our viscous and non inertial oscillator are included in \eqref{equadinn} where the $x-$equation, a Weber one, solved by means of the Hermite and Kummer functions, see
 \eqref{eh}. Our model has been analyzed carrying out five sample transients. In all cases the effects of the $\omega$ outside law, starting at $t=0$, are stopped at $t=10 s$. 
\figurename~\ref{f02} describes the $x$-oscillations whenever the angular speed $\omega$ is falling down: some viscous drag values have been considered in presence of a meaningful elastic force. The same problem has then been seen differently and we provide a polar plot of the bead torse's projection.
 The second transient shows the oscillating behavior of the constraint reaction $R_y$ under major elastic effects. The third case takes into account a linear $\omega$ growth: \figurename~\ref{f06} shows the oscillatory amplitudes' different increase according to different drag $A$ values. The fourth case displays a nonoscillatory behavior.The last one assumes the angular speed to be kept constant.
 Four of the above cases assume the rotating speed changing in time, different viscous drags, but no dry friction, namely $\mu=0$. Taking into account the dry friction on the bead too, the ODE keeps its linearity but looses its homogeneity. \figurename~\ref{finove} plots a set of damped oscillations with different $A$ values and forced by dry friction $\mu=1$. 
 
%\appendix 
\section*{Appendix}

The first hypergeometric series appeared in a Wallis's  book \cite{walli}: 
\[
_{2}{\rm F}_{1}\left( \left. 
\begin{array}{c}
a;b \\[2mm]
c
\end{array}
\right| x\right) =\sum_{n=0}^{\infty }\frac{\left( a\right) _{n}\left(
b\right) _{n}}{\left( c\right) _{n}}\frac{x^{n}}{n!} ,
\]
for $|x|<1$ and complex parameters $a,\,b,\,c$ being $(a)_n,\,(b)_n,\,(c)_n$ Pochhammer symbols, for instance:
\[
(a)_n=\frac{\Gamma(a+n)}{\Gamma(a)}
\] 
being $\Gamma$ the Euler factorial functions. Function $_{2}{\rm F}_{1}$ solves the linear second order Gauss differential equation for the unknown $u(x)$:
\begin{equation}\label{hy}\tag{A1}
x(1-x)\,u^{\prime\prime}+\left[c-\left(a+b+1\right)x\right]\,u^{\prime}-ab\,u=0
\end{equation}
Furthermore, many functions have been introduced in 19$^{\mathrm{th}}$ century either for
generalizing it to multiple variables or taking the special way of {\it confluence}. This is the case of the Kummer CHF, defined by the absolutely convergent infinite power series:
$$
_1{\rm F}_1\left( \left. 
\begin{array}{c}
a \\[2mm]
c
\end{array}
\right| x\right)=\sum_{n=0}^{\infty }\frac{\left( a\right) _{n}}{\left( c\right) _{n}}\frac{x^{n}}{n!}.
$$
It is analytic, regular at zero entire single-valued transcendental function of all $a,\,c,\,x,$ (real or complex) except $c = 0,\,-1,\, -2,\, -3,\,\ldots,$ for which it has simple poles. 
 The above series is a solution of the Kummer ordinary differential equation in $y(x)$:
\begin{equation}\label{Kumm}\tag{K}
xy'' +(c-x)y' -ay = 0,
\end{equation}
A detailed but compact outline of CHF integral representation, asymptotic and formul\ae\, of various functions through it, can be found in \cite{lebe}.

The French mathematician C. Hermite (1822-1901) considered the second order differential equation in $v(z)$:
\begin{equation}\label{hermy}\tag{A2}
v'' -2zv' +2\nu v = 0,
\end{equation}
The general solution of \eqref{hermy} is
\begin{equation}\label{hermysol}\tag{A3}
v(z)=c_1 H_{\nu }(z)+c_2 \,
   _1{\rm F}_1\left( \left. 
\begin{array}{c}
-\frac{\nu
   }{2} \\[2mm]
\phantom{-}\frac{1}{2}
\end{array}
\right| z^2\right)
\end{equation}
where $H_\nu (z)$ is the {\it Hermite function} of degree $\nu$ of the real variable $z$. If $\nu$ is an integer,  $H_\nu (z)$ reduces to the {\it Hermite polynomials} of degree $\nu$. Such polynomials are a certain subclass of the so called orthogonal polynomials and are met solving the simple harmonic oscillator of quantum mechanics. On the contrary, if $\nu$ is not an integer, as in our case, the Hermite function is a linear combination of Kummer functions, namely:
\[
H_\nu (z)=\sqrt{\pi }\, 2^{\nu } \left(\frac{1}{\Gamma \left(\frac{1-\nu
   }{2}\right)}\,\,
   _1{\rm F}_1\left( \left. 
\begin{array}{c}
-\frac{\nu
   }{2} \\[2mm]
\phantom{-}\frac{1}{2}
\end{array}
\right| z^2\right) -\frac{1}{\Gamma \left(-\frac{\nu
   }{2}\right)}2 z \,
   _1{\rm F}_1\left( \left. 
\begin{array}{c}
\frac{1-\nu}{2} \\[2mm]
\frac{3}{2}
\end{array}
\right| z^2\right)\right).
\]
Finally in the paper appeared a $_1{\rm F}_2$ function, whose power series is
\begin{equation}\label{1F2}\tag{A3}
_1{\rm F}_2\left( \left. 
\begin{array}{c}
a \\[2mm]
b_1,\,b_2
\end{array}
\right| x\right) =\sum_{n=0}^\infty\frac{(a)_n}{(b_1)_n(b_2)_n}\frac{x^n}{n!}
\end{equation}
which converges for any $x\in\mathbb{C}.$ Even if this function we did not appear in any differential in the paper, we point out that it is solution of the third order differential equation in the unknown $w(z)$
\[
z^2 w^{(3)}+\left(b_1+b_2+1\right) z w''+\left(b_1 b_2-z\right) w'-a_1 w=0.
\]

\subsection*{Acknowledgements}

\noindent The authors take the opportunity for thanking the referees for their constructive criticism.

\noindent The authors are indebted to professor Aldo Scimone who drew the first figure of this paper and warmly thank him. 

\noindent The last author is supported by an RFO grant issued by the Italian Ministry of University and research.

%\bibliographystyle{elsarticle-num}
%\bibliography{copter}

\end{document}